\documentclass[a4paper, notitlepage, 12pt]{article}

\usepackage[total={160mm,250mm}, top=25mm, left=25mm, includefoot]{geometry}

\usepackage{ucs}              
\usepackage[utf8x]{inputenc} 	

\usepackage[T1]{fontenc}     
\usepackage{lmodern}         

\usepackage[singlespacing]{setspace}

\usepackage{color}
\usepackage{amsmath,amsthm,amssymb,amscd,amsbsy}
\usepackage{upgreek}

\usepackage{tensor}

\usepackage{wasysym}
\usepackage{float}
\usepackage{tabularx,paralist}
\usepackage{graphicx}
\usepackage[breaklinks]{hyperref}
\usepackage{authblk}

\usepackage{multirow} 
\usepackage{array}
\usepackage{caption}

\usepackage{rcs} \RCS $Revision: 1.1.1.2 $ \RCS $Date: 2018/11/22 16:47:45 $ \RCS $Author: hgm $ \RCS $RCSfile: 18_Param-Lin_coupled.tex,v $ \RCS $Id: 18_Param-Lin_coupled.tex,v 1.1.1.2 2018/11/22 16:47:45 hgm Exp $

\usepackage[british]{babel}

\usepackage{refdef}
\usepackage{ifontdef}

\newcommand{\ignore}[1]{}

\newcommand{\Lp}{\mrm{L}}

\newcommand{\Hp}{\mrm{H}}
\newcommand{\Ck}{\mrm{C}}

\newcommand{\vrho}{\varrho}

\newcommand{\vsigma}{\varsigma}

\newcommand{\vphi}{\varphi}
\newcommand{\vpi}{\varpi}
\newcommand{\vkappa}{\varkappa}

\newcommand{\vek}[1]{\mathchoice{\displaystyle\boldsymbol{#1}}
{\textstyle\boldsymbol{#1}}{\scriptstyle\boldsymbol{#1}}
{\scriptscriptstyle\boldsymbol{#1}}}

\newcommand{\mat}[1]{\mathchoice{\displaystyle\mathbf{#1}}
{\textstyle\mathbf{#1}}{\scriptstyle\mathbf{#1}}
{\scriptscriptstyle\mathbf{#1}}}

\newcommand{\opb}[1]{\vek{{\mathsf{#1}}}}

\newcommand{\ops}[1]{\mathchoice{\displaystyle\mathsf{#1}}
{\textstyle\mathsf{#1}}{\scriptstyle\mathsf{#1}}
{\scriptscriptstyle\mathsf{#1}}}

\newcommand{\tnb}[1]{\mathchoice{\displaystyle\mathboldsans{#1}}
{\textstyle\mathboldsans{#1}}{\scriptstyle\mathboldsans{#1}}
{\scriptscriptstyle\mathboldsans{#1}}}

\newcommand{\tns}[1]{\mathchoice{\displaystyle\mathsans{#1}}
{\textstyle\mathsans{#1}}{\scriptstyle\mathsans{#1}}
{\scriptscriptstyle\mathsans{#1}}}

\newcommand{\diag}{\mathop{\mathrm{diag}}\nolimits}

\newcommand{\divg}{\mathop{\mathrm{div}}\nolimits}

\newcommand{\im}{\mathop{\mathrm{im}}\nolimits}
\newcommand{\dom}{\mathop{\mathrm{dom}}\nolimits}

\newcommand{\tr}{\mathop{\mathrm{tr}}\nolimits}

\newcommand{\spn}{\mathop{\mathrm{span}}\nolimits}

\newcommand{\dd}{\partial}
\newcommand{\di}{\mathrm{d}}

\newcommand{\KL}{Karhunen-Lo\`eve}
\newcommand{\ip}[2]{\langle #1, #2 \rangle}
\newcommand{\bkt}[2]{\langle #1 | #2 \rangle}

\newcommand{\nd}[1]{\| #1 \|}

\newcommand{\trpos}{{\ops{T}}}

\definecolor{myred}{rgb}{1, 0.2, 0.2}

\newcommand{\autheadcr}{\authorcr}
\newcommand{\citep}[1]{\cite{#1}}

\newcommand{\authorhgm}{Hermann G. Matthies}
\newcommand{\authorro}{Roger Ohayon}

\newcommand{\affilwire}{Institute of Scientific Computing\autheadcr
                        Technische Universit\"at Braunschweig\autheadcr
                        38092 Braunschweig, Germany\autheadcr
                        e-mail: \ttt{wire@tu-bs.de}}
\newcommand{\affilcnam}{Laboratoire de M\'ecanique des Structures et des
                       Syst\`emes Coupl\'es\autheadcr
                       Conservatoire National des Arts et M\'etiers (CNAM)\autheadcr
                      75141 Paris Cedex 03, France}

\newcommand{\thetitle}{Analysis of parametric models for coupled systems}

\newcommand{\theauthor}{\authorhgm, \authorro}
\newcommand{\thesubject}{35B30, 37M99, 41A05, 41A45, 41A63, 60G20, 60G60, 65J99, 93A30}
\newcommand{\thekeywords}{parametric models, reproducing kernel Hilbert space,
                          correlation, factorisation, spectral decomposition, representation,
                          coupled systems}

\newcommand{\textdate}{\today}

\newcommand{\thebib}{./bib}

\begin{document}

\title{\thetitle\thanks{Partly supported by the Deutsche
          Forschungsgemeinschaft (DFG) through SPP 1886 and SFB 880.}}


\author[*]{\authorhgm}
\author[o]{\authorro}


\affil[*]{\affilwire}
\affil[o]{\affilcnam}

\date{\textdate}


\maketitle

%

\begin{abstract}
In many instances one has to deal with parametric models.
Such models in vector spaces are connected to a linear map.
The reproducing kernel Hilbert space and affine- / linear- representations in terms
of tensor products are directly related to this linear operator.  This
linear map leads to a generalised correlation operator, in fact it provides
a factorisation of the correlation operator and of the reproducing kernel.

The spectral decomposition of the correlation and kernel, as well as the
associated \KL{}- or proper orthogonal decomposition are a direct consequence.
This formulation thus unifies many such constructions under a functional analytic
view.  Recursively applying factorisations in higher order tensor representations
leads to hierarchical tensor decompositions.

This format also allows refinements for cases when the parametric model has more
structure.  Examples are shown for vector- and tensor-fields with certain required
properties.  Another kind of structure is the parametric model of a coupled system.
It is shown that this can also be reflected in the theoretical framework.

\vspace{5mm}

\vspace{5mm}
{\noindent\textbf{AMS Classification:} \thesubject}

\end{abstract}

%
%
%
%
%









%

\section{Introduction}  \label{S:intro}
Parametric models are used in many areas of science, engineering, and
economics to describe variations or changes of some system.  They can
have many different uses such as evaluating the \emph{design} of some system,
or to \emph{control} its behaviour, or to \emph{optimise} the performance
in some way.  Another important case is when some of the parameters may
be uncertain and are modelled by random variables (RVs), and one
wants to perform \emph{uncertainty quantification} or \emph{identify}
some of the parameters in a mathematical model.  One important consideration
is the preservation of structure which one knows to be present in the
system.  One such structure is the
consideration of \emph{coupled systems}, and it will be shown how
this can be dealt with.  In fact, the coupling conditions can be
one of the possible parameters.

The representations of such parametric models leads directly to
\emph{reduced order models} which are used to lessen the possibly
high computational demand in some of the tasks described above.
Such reduced models hence become parametrised.
The survey \citep{BennWilcox-paramROM2015} and
the recent collection \citep{MoRePaS2015}, as well as
the references therein, provide a good account of parametric
reduced order models and some of the areas where they appear.
The interested reader may find there further information on
parametrised reduced order models and how to generate them.

This present work is a continuation of \citep{boulder:2011} and \citep{hgmRO-1-2018},
where the theoretical background of such parametrised models was treated
in a functional analysis setting.  For many of the theoretical details
we thus refer to these publications, and especially to \citep{hgmRO-1-2018}
for a more thorough account of the theory.

As an example, assume
that some physical system is investigated, which
is modelled by an evolution equation for its state $u(t) \in \C{V}$
at time $t \in [0,T]$,
where  $\C{V}$ is assumed to be a Hilbert space for the sake of
simplicity: $\dot{u}(t) = A(\mu;u(t)) + f(\mu;t);\quad u(0) = u_0$,
where the superimposed dot signifies the time derivative,
$A$ is an operator modelling the
physics of the system, and $f$ is some external excitation.
The model depends on some quantity $\mu\in\C{M}$, where $\C{M}$
denotes the set of possible parameters, and we assume that for
all $\mu$ of interest the system is well-posed.  Other than that, it
is not assumed that the set $\C{M}$ has any additional structure.
In this way the system state becomes a function of the parameters,
and can thus be written as $u(\mu;t)$.  Later mainly the dependence on
$\mu$ will be interesting, so that the other arguments may be dropped.

To fix ideas, take as a simple example a very simple fluid-structure
interaction problem used in \citep{LefrancoisBoufflet2010} to explain
the basics of coupling algorithms: it is a mass-spring system coupled
with a gas-filled piston.  The governing equations
are in the simplest case considered there (with slight change of notation):
\[ 
m \,\ddot{w}(t) + k\, w(t)  =  S (p(t) - p_0); \qquad
p(t) = p_0 \left( 1+\frac{\gamma-1}{2} \frac{\dot{w}(t)}{c_0}\right)^{\frac{2\gamma}{\gamma-1}}.
\]
The mass $m$ and spring constant $k$ are properties of the mass-spring system
with displacement $w(t)$
coupled over the surface $S$ with the gas-filled piston with pressure $p(t)$.
The equilibrium pressure is $p_0$, and $c_0$ is the speed of sound in the gas
with specific heat ratios $\gamma$.  Introducing the velocity $v(t)=\dot{w}(t)$,
the state of the system is given by $u = [w, v]^\trpos$, the system depends on the
parameters $\mu=(m,k,S,c_0,\gamma-1)$, the state Hilbert space is $\C{V}=\D{R}^2$,
the set $\C{M}$ can be taken as $\C{M}=\{ (\mu_j)_{j=1..5} \mid \mu_j > 0 \}
\subset \D{R}^5$, and the problem
$\dot{u}(t) = A(\mu;u(t)) + f(\mu;t)$ becomes:
\begin{multline} \label{eq:mkp}
\dot{u}(t) = \begin{bmatrix}\dot{w}(t) \\ \dot{v}(t) \end{bmatrix} = 
\begin{bmatrix} v(t) \\ -\frac{k}{m} w(t) + \frac{S p_0}{m} 
\left( 1+\frac{\gamma-1}{2} \frac{v(t)}{c_0}\right)^{\frac{2\gamma}{\gamma-1}}\end{bmatrix}
+ \begin{bmatrix} 0 \\ - \frac{S p_0}{m}\end{bmatrix}.
\end{multline}
A little twist can be given to this by assuming that some, or all, of the
parameters are \emph{random variables}; say for example the spring stiffness $k$.
Formally, this is a measurable function from some probability space into
the real numbers: $k:\Omega\to\D{R}$.  This makes also the displacement $w$ and
velocity $v$ into random variables $w, v: \Omega\to\D{R}$.  If we assume that all
involved random variables have finite variance, i.e.\ $k, w, v \in \Lp_2(\Omega)$,
then the parameter set could be taken as $\C{M} = \D{R}_+ \times \Lp_2(\Omega) \times\D{R}_+^3$,
and the state space would be the Hilbert space $\C{V}=\Lp_2(\Omega)^2$.  Such
probabilistic examples have prompted much of the theory and terminology \citep{kreeSoize86},
and such probabilistic problems are treated specifically in the present framework
in \citep{hgm-3-2018}; but here we do not want to digress and keep the focus on coupled systems.

A bit more involved is the following example from \citep{MatthiesSteindorf2003},
it is a kind of generic example of fluid-structure interaction.  The fluid is described by the
incompressible Navier-Stokes equation in arbitrary Lagrangean-Eulerian (ALE)
formulation in a domain $\Omega_f$:
\begin{align} \label{eq:FSI-fl-1}
\rho_f(\dot{v} + ((v-\dot{\chi})\cdot\nabla) v) - \divg \sigma + \nabla p &= r_f, \\
\label{eq:FSI-fl-2} 2 \sigma = \nu_f (\nabla v + (\nabla v)^{\trpos}),\quad \divg v &= 0,
\end{align}
where $\rho_f$ and $\nu_f$ are the fluid mass-density and viscosity, $v(x,t)$ is the fluid
velocity-field, $\sigma(x,t)$ is the viscous stress in the fluid, $r_f(x,t)$ are the
volume forces in the fluid, and the pressure
$p(x,t)$ is the Lagrange multiplier for the incompressibility constraint, whereas
$\chi(x,t)$ is the movement of the ALE background reference system.  On part of
the boundary $\Gamma_c\subset\dd\Omega_f$ the fluid is coupled to an elastic solid,
described in the solid domain $\Omega_s$ (with also $\Gamma_c\subset\dd\Omega_s$) in
a Lagrangean or material frame by
\begin{align} \label{eq:FSI-sd-1}
\rho_s\ddot{w}  - \divg (\tns{F} \tns{S}) = r_s, \quad \tns{F} &= \nabla w, \\
\label{eq:FSI-sd-2} \tns{S} = \lambda_s (\tr \tns{E}) \tns{I} + 
  2\mu_s \tns{E},\quad 2\tns{E}=\tns{C}-\tns{I}, \quad \tns{C}&= \tns{F}^{\trpos}\tns{F},
\end{align}
where  $\rho_s$ is the solid mass-density and $\lambda_s, \mu_s$ are the elastic Lamé
moduli, $w(x,t)$ is the solid displacement-field and $\tns{F}$ its gradient,
$\tns{S}(x,t)$ is the elastic 2nd Piola-Kirchhoff stress in the solid, $r_s(x,t)$ are the
solid volume forces, and the Lagrange-Green strain $\tns{E}$ is stated in terms of the Cauchy-Green
strain tensor $\tns{C}$.  The position of a solid particle which was at position $x$ at
time $t=0$ is $w(x,t) + x$.  Hence at the each point $\chi(x,t)=w(x,t) + x$
on the coupling boundary $\Gamma_c$  with normal $\vek{n}(x,t)$ one has
the condition that the velocities of fluid and solid have to match
\begin{equation} \label{eq:FSI-cpl}
  v(\chi(x,t),t) = \dot{w}(x,t).
\end{equation}
The solution $(v,p)$ to the fluid part \refeq{eq:FSI-fl-1} and \refeq{eq:FSI-fl-2}
lives in a Hilbert space $\Hp(\divg)(\Omega_f)\times\Lp_2(\Omega_f)$, whereas the
displacement $w$ of the solid part can be envisioned in $\Hp^1(\Omega_s)^3$ with a velocity
$\dot{w}\in\Hp^1(\Omega_s)^3$.  The state of the system $(v,p,w,\dot{w})$
is thus in the Hilbert space $\C{V} = \Hp(\divg)(\Omega_f)\times\Lp_2(\Omega_f)
\times \Hp^1(\Omega_s)^6$.  The parameters can for example be the material
constants for fluid and solid $\mu = (\rho_f,\nu_f,\rho_s,\lambda_s,\mu_s)$ such that
the set $\C{M}$ can be taken again as $\C{M}=\{ (\mu_j)_{j=1..5} \mid \mu_j > 0 \}
\subset \D{R}^5$, or it can be the smooth initial shape $\chi_0(x)$ of the coupling boundary
$\Gamma_c$, such that $\mu = \chi_0$ and $\C{M}\subset \Ck^1(\Omega_c)$ with the
boundary description $\Omega_c\subset\D{R}^2$; or a combination of these two cases.

With these two examples in mind, turning again to the general description,
one is interested in how the system changes when these
\emph{parameters} $\mu$ change.  As we have seen in the above examples, these
can be something specific describing the operator, or the initial condition, or
specifying the excitation, etc.\ \citep{BennWilcox-paramROM2015}.
One may be interested in the state of the system $u(\mu;t)$,
or some functional of it, say $\Psi(\mu)$.  In the first example this could
for example be the maximum acceleration $\Psi(\mu)=\max_t \ddot{w}(\mu;t)$.
While evaluating $A(\mu;\hat{u}(t))$ --- for some \emph{trial} state $\hat{u}(t)$ ---
or $f(\mu;t)$ for a certain $\mu$ may be straightforward,
there are situations where
evaluating $u(\mu;t)$ or $\Psi(\mu)$ may be very costly.

In this situation one is interested in representations of $u(\mu;t)$ or $\Psi(\mu)$
which allow a cheaper evaluation, these
are called \emph{proxy}- or \emph{surrogate}-models, among others.
Any such parametric object can be analysed by linear maps which are associated with
such representations.
This association of parametric models and linear mappings has probably been
known for a long time, see \citep{kreeSoize86} for an exposition
in the context of stochastic models.

It should be pointed out that the exposition here is a general
technique which can be used to analyse any parametric model and its approximations,
where the parameters can be of whatever nature, e.g.\ just some numbers,
or functions, or random variables, etc.,
and not so much something to be directly implemented.  It is shown though
that in the framework of orthogonal bases it has a direct connection
with proper orthogonal decomposition (POD) and \KL{} expansions. 
As will be seen, it also connects with representations in tensor
product spaces, which allows numerically to use low-rank tensor approximations
\citep{Hackbusch_tensor, boulder:2011}.  It is furthermore also connected
with non-orthogonal decompositions which are easier to compute,
like the \emph{proper generalised decomposition}(PGD) \citep{chinestaBook}.
Here we shall only consider orthogonal bases for the sake of conciseness
of exposition.  

Whereas the parametric map may be quite complicated, the association
with a linear map translates the whole problem into one of linear
functional analysis, and into linear algebra upon approximation and
actual numerical computation.

%
%
%
%
%
%
%
%
%
%
%



\section{Parametric models and reproducing kernel} \label{S:parametric}
This is a short recap of the developments in \citep{hgmRO-1-2018}, where
the interested reader may find more detail.
Let $r: \C{M} \rightarrow \C{U}$ be
one of the objects alluded to in the introduction, where $\C{M}$ is
some set without any further assumed structure,
and $\C{U}$ is assumed for the sake of simplicity
as a separable Hilbert space with inner product
$\langle\cdot|\cdot\rangle_{\C{U}}$.

Assume without significant loss of generality that
$\spn r(\C{M}) = \spn \im r \subseteq \C{U} $,
the subspace of $\C{U}$ which is spanned by all the
vectors $\{r(\mu)\mid  \mu\in\C{M}\}$,  is dense in $\C{U}$ ---
otherwise we just restrict ourselves to the closure of $\spn r(\C{M}) = \spn \im r$.
Then to each such function $r$ one may associate a linear map
$\tilde{R}: \C{U} \ni u \mapsto \bkt{r(\cdot)}{u}_{\C{U}} \in \D{R}^{\C{M}}$, the
real-valued functions on the set $\C{M}$.  As a motivation,
one may think of such functions as providing a `co-ordinate system'
on the otherwise unstructured set $\C{M}$.
By construction, $\tilde{R}$ restricted to $\spn\im r = \spn r(\C{M})$ is injective,
and has an inverse on its restricted range range
$\tilde{\C{R}}:= \tilde{R}(\spn \im r) \subseteq \D{R}^{\C{M}}$.
This may be used to define an inner product on $\tilde{\C{R}}$ as
\begin{equation}     \label{eq:V}
\forall \phi, \psi \in \tilde{\C{R}} \quad \langle \phi | \psi \rangle_{\C{R}} :=
     \langle \tilde{R}^{-1} \phi | \tilde{R}^{-1} \psi \rangle_{\C{U}},
\end{equation}
and to denote the completion of $\tilde{\C{R}}$ with this inner product by $\C{R}$.
One immediately obtains that $\tilde{R}^{-1}$ is a bijective isometry between $\spn \im r$
and $\tilde{\C{R}}$, hence extends to a \emph{unitary} map between $\C{U}$ and $\C{R}$,
as does $\tilde{R}$.

Given the maps $r:\C{M}\to\C{U}$ and $\tilde{R}:\C{U}\to\C{R}$, one may define the
\emph{reproducing kernel} \citep{berlinet, Janson1997}
given by $\vkappa(\mu_1, \mu_2) := \bkt{r(\mu_1)}{r(\mu_2)}_{\C{U}}$.
It is straightforward to verify that $\vkappa(\mu,\cdot)\in\tilde{\C{R}}\subseteq\C{R}$,
and $\spn \{ \vkappa(\mu,\cdot)\;\mid\; \mu\in\C{M} \}=\tilde{\C{R}}$,
as well as the reproducing property
$\phi(\mu) = \bkt{\vkappa(\mu,\cdot)}{\phi}_{\C{R}}$ for all $\phi\in\tilde{\C{R}}$.

On this \emph{reproducing kernel Hilbert space} (RKHS) $\C{R}$
one can build a first representation.
As $\C{U}$ is separable, so is $\C{R}$, 
and one may choose a complete orthonormal system (CONS) $\{\vphi_m\}_{m\in\D{N}}$
in $\C{R}$.  Then with the CONS
$\{y_m \;|\; y_m= \tilde{R}^{-1} \vphi_m = \tilde{R}^*\vphi_m\}_{m\in\D{N}}$ in $\C{U}$,
the unitary operator $\tilde{R}$, and its adjoint or inverse $\tilde{R}^*=\tilde{R}^{-1}$,
and the parametric element $r(\mu)$ become \citep{hgmRO-1-2018}
\begin{multline} \label{eq:VII0}
  \tilde{R} = \sum_m \vphi_m \otimes y_m;  \quad \text{i.e. } \quad 
  \tilde{R}(u)(\cdot) = \sum_m \bkt{y_m}{u}_{\C{U}} \vphi_m(\cdot), \quad
  \tilde{R}^*=\tilde{R}^{-1} = \sum_m y_m \otimes \vphi_m;\\
   r(\mu) = \sum_m \vphi_m(\mu) y_m = \sum_m \vphi_m(\mu)\, \tilde{R}^* \vphi_m .
\end{multline}
Observe that the relations \refeq{eq:VII0} exhibit the tensorial nature of the
representation mapping.
One sees that \emph{model reductions} may be achieved by choosing only subspaces of $\C{R}$,
i.e.\ a---typically finite---subset of $\{\vphi_m\}_{m}$.
Furthermore, the representation of $r(\mu)$ in \refeq{eq:VII0} is \emph{linear}
in the new `parameters' $\vphi_m$.

%
%
%
%
%
%
%
%


%

\section{Correlation and kernel space}  \label{S:correlat}
The RKHS construction $\C{R}$ of \refS{parametric}
just mirrors or reproduces the 
inner product structure on the original space $\C{U}$.
There is presently no way of telling what is important in the parameter set $\C{M}$.
For this one needs additional information.
As a way of indicating what is important on the set $\C{M}$, assume that
there is another inner product $\bkt{\cdot}{\cdot}_{\C{Q}}$ for
scalar functions $\phi \in \D{R}^{\C{M}}$, and denote the Hilbert
space of functions with that inner product by $\C{Q}$.  We define a
linear map in the same way as the map $\tilde{R}$ in \refS{parametric},
but as it now has a different range or image space --- $\C{Q}$ instead of
$\C{R}$ --- and denote the map by $R:\C{U}\ni u\mapsto \bkt{r(\cdot)}{u}_{\C{U}}\in\C{Q}$.
As the inner product on the image space has changed, the map $R$, unlike the map $\tilde{R}$
in  \refS{parametric}, will in general not be unitary any more.
Also assume that the subspace $\dom R = \{ u\in\C{U} \mid \nd{Ru}_{\C{Q}}< \infty\}$
is, if not the whole space $\C{U}$, at least dense in $\C{U}$, and that
the densely defined operator $R$ is closed.  These are essentially requirements
that the topologies on $\C{U}$ and
$\C{Q}$ fit in some way with the map $r:\C{M}\to\C{U}$.
But to make things even simpler, assume here that $R$ is
defined on the whole space and hence continuous, and still injective,
unless explicitly otherwise stated.

With this, one may define \citep{kreeSoize86, hgmRO-1-2018}
a densely defined map $C$ in $\C{U}$ through the bilinear form 
\begin{equation}   \label{eq:IX}
   \forall u, v \in \C{U}:\quad \bkt{Cu}{v}_{\C{U}} := \bkt{Ru}{Rv}_{\C{Q}} .
\end{equation}
The map $C$, which may also be written as $C=R^* R$, may be called the 
\emph{`correlation'} operator.  By construction it is self-adjoint
and positive, and if $R$ is continuous so is $C$.
In case the inner product $\bkt{\cdot}{\cdot}_{\C{Q}}$ comes from a
measure $\vpi$ on $\C{M}$, so that for two functions
$\phi$ and $\psi$ on $\C{M}$ one has
\[
  \bkt{\phi}{\psi}_{\C{Q}} := \int_{\C{M}} \phi(\mu) \psi(\mu) \; \vpi(\di \mu),
  \quad \text{such that from \refeq{eq:IX}}
\]
\[  
  \bkt{Cu}{v}_{\C{U}} = 
  \int_{\C{M}} \bkt{r(\mu)}{u}_{\C{U}}\bkt{r(\mu)}{v}_{\C{U}} \; \vpi(\di \mu),\text{ i.e. }
  C = R^* R = \int_{\C{M}} r(\mu) \otimes r(\mu) \; \vpi(\di \mu).
\]
The space $\C{Q}$ may then be taken as $\C{Q}:=\Lp_2(\C{M},\vpi)$.
A special case is when $\vpi$ is a probability measure,
$\vpi(\C{M}) = 1$, this inspired the term `correlation'  \citep{kreeSoize86}.

In \refS{parametric} it was the factorisation of $C= R^*R$ which allowed
the RKHS representation in \refeq{eq:VII0}.   For other representations, one needs other
factorisations.  Most common is to use the spectral
decomposition (e.g.~\citep{DautrayLions3}) of $C$ to achieve such a factorisation.

On infinite dimensional Hilbert spaces self-adjoint operators may have
a continuous spectrum, e.g.~\citep{DautrayLions3}.  To make everything as simple as possible
to explain the main underlying idea,
we shall from now on assume that $C$ is a non-singular trace class or nuclear operator.
This means that it is compact, the spectrum $\sigma(C)$ is a point spectrum, has a
CONS $\{v_m\}_m$ consisting of eigenvectors, with
each eigenvalue $\lambda_m\ge \lambda_{m+1}\dots\ge 0$ positive
and counted decreasingly according to
their finite multiplicity, and has finite trace $\tr C = \sum_m \lambda_m < \infty$.
Then a version of the spectral decomposition of $C$ is
\begin{equation}   \label{eq:XIII}
C  = \sum_m \lambda_m (v_m \otimes v_m) .
\end{equation}
Define a new CONS $\{s_m\}_m$ in $\C{Q}$:
 $\lambda_m^{1/2} s_m := R v_m$, to obtain the corresponding
\emph{singular value decomposition} (SVD) of $R$ and $R^*$:
\begin{multline}   \label{eq:XIV}
R = \sum_m \lambda_m^\frac{1}{2} (s_m \otimes v_m)\,;   \quad \text{i.e. } \quad 
  R(u)(\cdot) = \sum_m \bkt{v_m}{u}_{\C{U}} s_m(\cdot), \quad
R^* = \sum_m \lambda_m^\frac{1}{2} (v_m \otimes s_m)\,; \\
r(\mu) =  \sum_m \lambda_m^\frac{1}{2} \, s_m(\mu) v_m = \sum_m s_m(\mu)\, R^* s_m,
\text{ as } \; R^*s_m = \lambda_m^\frac{1}{2} \, v_m .
\end{multline}
The set $\vsigma(R)=\{\lambda_m^{1/2}\}_m = \sqrt{\sigma(C)}\subset \D{R}_+$
are the \emph{singular values} of $R$ and $R^*$.
The last relation is
the so-called \emph{\KL{} expansion} or \emph{proper orthogonal decomposition} (POD).
The finite trace condition of $C$ translates into the fact
that $r$ is in $\C{U}\otimes\C{Q}$.
If in that relation the sum is \emph{truncated} at $n\in\D{N}$, i.e.\
\begin{equation}   \label{eq:best-n-term}
r(\mu) \approx r_n(\mu) =  \sum_{m=1}^n \lambda_m^\frac{1}{2} \, 
        s_m(\mu) v_m = \sum_{m=1}^n s_m(\mu)\, R^* s_m,
\end{equation}
we obtain the \emph{best $n$-term approximation} to $r(\mu)$ in the norm
of $\C{U}$.

Observe that, similarly to \refeq{eq:VII0}, $r$ is linear in the $s_m$.
This means that by choosing the `co-ordinate transformation'
$\C{M}\ni \mu\mapsto (s_1(\mu),\dots,s_m(\mu),\dots)\in\D{R}^{\D{N}}$ one obtains
a \emph{linear / affine} representation where the first co-ordinates are
the most important ones.

A formulation of the spectral decomposition different from \refeq{eq:XIII} does
not require $C$ to be nuclear \citep{DautrayLions3},
nor does $C$ or $R$ have to be continuous.
The self-adjoint and positive operator $C:\C{U}\to\C{U}$ is unitarily equivalent
with a multiplication operator $M_{\gamma}$,
\begin{equation}  \label{eq:ev-mult}
  C = V M_{\gamma} V^*,
\end{equation}
where $V:\Lp_2(\C{T})\to\C{U}$ is unitary between some $\Lp_2(\C{T})$ on a measure
space $\C{T}$ and the Hilbert space $\C{U}$, and $M_{\gamma}$ is a multiplication operator,
multiplying any $\psi\in\Lp_2(\C{T})$ with a real-valued
function $\gamma$.
In case $C$ is bounded, so is $\gamma\in\Lp_\infty(\C{T})$.
As $C$ is positive, $\gamma(t)\ge 0$ for $t\in\C{T}$,
the essential range of $\gamma$ is the spectrum of $C$, and
$M_{\gamma}$ is self-adjoint and positive.  Its square root is
$M_{\gamma}^{1/2} := M_{\sqrt{\gamma}}$, from which one obtains
the square-root of $C^{1/2} = V  M_{\sqrt{\gamma}} V^*$.
The factorisation corresponding to $C=R^* R$ with the square-root is
$C= (V M_{\sqrt{\gamma}})(V M_{\sqrt{\gamma}})^* =: G^* G$.
Another possibility is $C = (C^{1/2})^* C^{1/2} = C^{1/2} C^{1/2}$.
From this follows another formulation of the singular value decomposition
(SVD) of $R$ and $R^*$ with a unitary $U:\Lp_2(\C{T})\to\C{Q}$:
\begin{equation}  \label{eq:ev-mult-svd}
R = U M_{\sqrt{\mu}} V^*,\quad R^* = V M_{\sqrt{\mu}} U^*.
\end{equation}

These are all examples of a general factorisation
$C = B^* B$, where $B:\C{U}\to\C{H}$ is a map to a Hilbert space $\C{H}$
with all the properties demanded from $R$---see the beginning of this section.
It can be shown \citep{hgmRO-1-2018} that
any two such factorisations $B_1:\C{U}\to\C{H}_1$ and $B_2:\C{U}\to\C{H}_2$ 
with $C=B_1^*B_1=B_2^*B_2$ are
\emph{unitarily equivalent} in that there is a unitary map $X_{21}:\C{H}_1\to\C{H}_2$
such that $B_2 = X_{21} B_1$.  Equivalently, each such factorisation is unitarily
equivalent to $R$, i.e.\  there is a unitary $X:\C{H}\to\C{Q}$ such that
$R= X B$.
For finite dimensional spaces, a favourite choice 
is the Cholesky factorisation $C = L L^{\trpos}$, where $B=L^{\trpos}$ and $B^* = L$.

In the situation where $C$ has a purely discrete spectrum
and a CONS of eigenvectors $\{v_m\}_m$ in $\C{U}$,
the map $B$ from the decomposition $C=B^* B$ can be used to
define a CONS $\{h_m\}_m$ in $\C{H}$:  $h_m := B C^{-1/2} v_m$, which is
 an eigenvector CONS of the operator $C_{\C{H}} := B B^*:\C{H}\to\C{H}$,
with $C_{\C{H}} h_m := \lambda_m h_m$, see   \citep{hgmRO-1-2018}.
From this follows a SVD of $B$ and $B^*$ analogous to \refeq{eq:XIV}.
Taking the special case $\C{H}=\C{Q}$ with $C_{\C{Q}} = R R^*$, we see
that $C_{\C{Q}} s_m = \lambda_m s_m$, and $s_m = U V^* v_m$, as
well as $C_{\C{Q}} = U V^* C V U^*$.

From the factorisation and \KL{} expansion in \refeq{eq:XIV} one has
$r(\mu) = \sum_m s_m(\mu)\, R^* s_m$.  With other equivalent
factorisations $C=B^* B$ one obtains new representations in an
analogous manner.
The main result is  \citep{hgmRO-1-2018} that in the case of a nuclear
$C$ every factorisation leads to a separated representation in terms
of a series, and vice versa.
The associated `correlations' $C_{\C{Q}} = R R^*$ on $\C{Q}$ resp.\  $C_{\C{H}} = B B^*$
on $\C{H}$ have the same spectrum as $C$, and factorisations of $C_{\C{Q}}$ resp.\ $C_{\C{H}}$
induce new factorisations of $C$.

The abstract equation $C_{\C{Q}} = U V^* C V U^* = U M_\gamma U^* = 
\sum_m \lambda_m s_m \otimes s_m$ can be spelt out in more analytical
detail for the special case when the inner product on $\C{Q}$ is given
by a measure $\vpi$ on $\C{P}$.  It then becomes for all $\vphi, \psi\in\C{Q}$:
\begin{multline*}  
   \bkt{C_{\C{Q}}\vphi}{\psi}_{\C{Q}} = \sum_m \lambda_m \bkt{\vphi}{s_m}_{\C{Q}}
   \bkt{s_m}{\psi}_{\C{Q}} =    \bkt{R^* \vphi}{R^* \psi}_{\C{U}} =\\ 
   \iint_{\C{M}\times\C{M}}  \vphi(\mu_1) \vkappa(\mu_1, \mu_2)
    \psi(\mu_2)\; \vpi(\di \mu_1) \vpi(\di \mu_2) = \\
   \iint_{\C{M}\times\C{M}}  \vphi(\mu_1) \left( \sum_m \lambda_m s_m(\mu_1)  s_m(\mu_2)
    \right) \psi(\mu_2)\; \vpi(\di \mu_1) \vpi(\di \mu_2) = \\ 
     \sum_m \lambda_m \left(  \int_{\C{M}} \vphi(\mu) s_m(\mu)\; \vpi(\di \mu) \right) 
     \left( \int_{\C{M}} s_m(\mu) \psi(\mu) \vpi(\di \mu) \right),
\end{multline*}
i.e.\ $C_{\C{Q}}$ is a Fredholm integral operator and its spectral decomposition
is nothing but the familiar theorem of Mercer \citep{courant_hilbert} for the kernel
$\vkappa(\mu_1, \mu_2)=\sum_m \lambda_m s_m(\mu_1)  s_m(\mu_2)$.  Factorisations
of $C_{\C{Q}}$ are then usually factorisations of the kernel $\vkappa(\mu_1, \mu_2)$.

An example is if on some measure space
$(\C{X},\nu)$ it holds that $\vkappa(\mu_1, \mu_2)
 = $ \\ $\int_{\C{X}} g(\mu_1,x) g(\mu_2,x)\; \nu(\di x)$,
then the integral transform with kernel $g$ will play the role of a
factor as before did the mappings $R$ or $B$, leading to a new
\KL-like representation of $r$. 
The abstract setting outlined in this section can hence be applied to the analysis
of a great number of different situations, of any kind of representation of
$r(\mu)$, see \citep{hgmRO-1-2018} for more detail, and \citep{hgm-3-2018} for
the case where essentially $\mu$ is a random variable.

Let us remark how this framework here may be extended when the
assumptions at the beginning of this section are not satisfied.  In case 
the parametric map $r:\C{M}\to\C{U}$ maps not into a Hilbert space,
but only into a locally convex topological vector space (LCTVS), it is still possible to
define a corresponding linear map and a correlation operator.  
Denote the continuous dual space
of $\C{U}$ by $\C{U}^*$ and define a linear map by $S:\C{U}^* \ni u^* \mapsto S(u^*) :=
\ip{u^*}{r(\cdot)}_* \in \C{Q}$, where $\ip{\cdot}{\cdot}_*$ is the duality pairing on
$\C{U}^*\times\C{U}$.  The Hilbert space $\C{Q}$ is still identified with
its dual, although one could also make generalisations here.  Then the dual
map $S^*:\C{Q}\to\C{U}$ --- w.r.t.\ the weak* topology on $\C{U}^*$ ---
is defined as usual by $\ip{u^*}{S^*\phi}_* :=\bkt{Su^*}{\phi}_{\C{Q}}$.
This allows us to define a `correlation' $C_S:\C{U}^*\to\C{U}$ by $C_S = S^* S$,
i.e.\ $\ip{u^*}{C_S v^*}_* := \bkt{S u^*}{S v^*}_{\C{Q}}$.  The dual map $S^*:\C{Q}\to\C{U}$
can still provide a representation.  Other factorisations of $C_S$ such as $C_S = B^* B$, where
$B:\C{U}^*\to\C{H}$ --- a Hilbert space identified with its dual --- can provide
alternative representations via the map $B^*:\C{H}\to\C{U}$.   But spectral theory
can not be used easily as we are not in a Hilbert space and domain and range or image space
are not the same for $C_S$.

A frequent situation where some further development is possible
is as follows: there is an injective continuous map
$T:\C{Z}\to\C{U}$ from a Hilbert space $\C{Z}$, which will be a
pivot space identified with its dual into the space $\C{U}$,
such that the subspace $T(\C{Z})\subseteq \C{U}$ is dense.  Then the dual map
$T^*:\C{U}^* \to \C{Z}$ is also injective and continuous, $\C{T} :=T^*(\C{U}^*)\subseteq\C{Z}$
is dense in $\C{Z}$, and we have a Gel'fand-like triplet of spaces
\[  \C{U}^* \overset{T^*}{\to} \C{Z} \overset{T}{\to} \C{U}. \]
On the dense subspace $\C{T}\subseteq\C{Z}$ the map $T^*$ is invertible,
$J := T_{|\C{T}}^{-*}:\C{T} \to \C{U}^*$.
This allows one to define a mapping $R$ and a `correlation' $C$ densely defined on
$\C{T}\subseteq\C{Z}$:
\[ 
 R:=S \circ J : \C{T} \overset{J}{\to} \C{U}^* \overset{S}{\to} \C{Q} \quad \text{and} \quad
 C:= R^* R = J^* S^* S J = J^* C_S J: \C{T}\to\C{Z}.
\]
We view $C$ now as a self-adjoint positive operator densely defined in the Hilbert space $\C{Z}$,
and we are in the previously described Hilbert space setting.  

A typical occurance of such a situation is the case when $\C{Z} \hookrightarrow \C{U}$
is a \emph{continuously embedded} Hilbert space, the map $T$ is then just the identity.
A concrete example of this is $r:\C{M}\to\E{S}^\prime(\Omega)$, a parametric map into the
Schwartz space of tempered distributions.  As the dual of $\C{U}:=\E{S}^\prime(\Omega)$
in the weak* topology is $\C{U}^* = \E{S}(\Omega)$, the test space of rapidly decaying
smooth functions, we have the continuous embeddings $\C{U}^*=\E{S}(\Omega)
\hookrightarrow \C{Z}:=\Lp_2(\Omega) \hookrightarrow \E{S}^\prime(\Omega) = \C{U}$. 
Hence one may define the densely
defined maps $R:\C{T}:=\E{S}(\Omega)\to \C{Q}$ and $C: \C{Z}:=\Lp_2(\Omega) \to \C{Z}$,
given for $f, g \in \E{S}(\Omega)$ by the bilinear form $\bkt{Cf}{g}_{\Lp_2} = 
\bkt{R f}{R g}_{\C{Q}}=  \bkt{\ip{r(\cdot)}{f}_{\E{S}}}{\ip{r(\cdot)}{g}_{\E{S}}}_{\C{Q}}$,
where $\ip{\cdot}{\cdot}_{\E{S}}$ is the duality pairing on $\E{S}^\prime(\Omega) \times
\E{S}(\Omega)$.

%
%
%
%
%
%
%
%
%
%
%
%
%


%

\section{Structure preservation and coupled systems}  \label{S:xmpls}
The main feature up to now was the mapping $r:\C{M}\to\C{U}$ and the associated
linear map $R:\C{U}\to\C{Q}\subseteq\D{R}^{\C{M}}$, and the resulting tensor representation
of $r \in \C{U}\otimes\C{Q}$.  Here we mention some possible refinements and extensions
which try to make some known structure in the parametric model explicitly visible.

One frequent situation is that the parameter space is a product space, say
$\C{M} = \C{M}_1 \times \C{M}_2$ and a corresponding factorisation of the space
$\C{Q} = \C{Q}_1 \otimes \C{Q}_2$, where $\C{Q}_1\subseteq\D{R}^{\C{M}_1}$ and
$\C{Q}_2\subseteq\D{R}^{\C{M}_2}$.  Then the functions $\phi(\mu)\in\C{Q}$ are
linear combinations of products $\vphi(\mu_1)\psi(\mu_2)$, with $\vphi(\mu_1)\in\C{Q}_1$
and $\psi(\mu_2)\in\C{Q}_2$.  Now such a product may be seen as a parametric mapping
\[
 \vrho:\C{M}_2\ni \mu_2 \mapsto \vphi(\cdot)\psi(\mu_2) \in \C{Q}_1;
\]
and by setting
\begin{equation}  \label{eq:XXIII-tens}
  \C{Q} = \C{U}_*\otimes \C{Q}_* = \C{Q}_1 \otimes \C{Q}_2,
\end{equation}
the theory of the preceding sections may now be applied to the parametric map
$\vrho$ and its tensor product representation in $\C{U}_*\otimes \C{Q}_*$ to get
a refined representation of the complete model.  This shows how the product structure of
the parameter set $\C{M} = \C{M}_1 \times \C{M}_2$ expresses itself in the representation.
In some way the full parameter set is a \emph{coupled} object, and some of this is
reflected in the factorisation.  But for a real coupled system we will demand
a bit more, as will be explained later.

It is now not difficult to see that in case $\C{M} = \prod_j \C{M}_j$ with corresponding
$\C{Q} = \bigotimes_j \C{Q}_j$ and $j\in\C{J}\subset\D{N}$,
this may be further factorised by different associations depending on a
partition of the parameter set $\C{J}=\C{J}_1 \cup \C{J}_2$ into two disjoint sets
$\C{J}_1 \cap \C{J}_2 =\emptyset$:
\begin{equation}  \label{eq:XXIII-tens-p}
  \C{Q} = \C{U}_{*}\otimes\C{Q}_{*} = \left(\bigotimes_{k\in\C{J}_1} \C{Q}_k\right)
   \otimes \left(\bigotimes_{k\in\C{J}_2} \C{Q}_k\right) .
\end{equation}
Each of the factors can then be recursively factorised further, and
this leads to hierarchical tensor approximations, e.g.\ 
\citep{Hackbusch_tensor, boulder:2011}.

Of course it is possible to split the tensor product in different ways, and
the grouping of indices can be viewed as a tree.  The well-known
\emph{canonical polyadic} (CP) decomposition 
uses the flat tensor product in \refeq{eq:XXIII-tens-p}.
It has also been published as so-called \emph{proper generalised decomposition}
(PGD) as a computational method to solve high-dimensional
problems, see the review
\citep{chinestaPL2011} and the monograph \citep{chinestaBook}.
But recursive splittings of \refeq{eq:XXIII-tens-p} yield \emph{deep} or
\emph{hierarchical} tensor approximations.
Particular formats are the \emph{tensor train} (TT) and more
generally the \emph{hierarchical Tucker} (HT) decompositions,
see the review \citep{GrasedyckKressnerTobler2013} and the
monograph \citep{Hackbusch_tensor}.
These hierarchical low-rank tensor representations are
connected with deep neural networks \citep{CohenSha2016, KhrulkovEtal2018}.
It is the eigenvalue structure of the correlation $C$ or equivalently the
structure of the singular values of the associated linear map $R$ in the particular
splitting \refeq{eq:XXIII-tens-p} which determines how many terms a series
representation needs to be a good reduced model with a certain accuracy.

Another frequent case is that the role of the Hilbert space $\C{Q}$ is taken by a
tensor product $\C{W}=\C{Q} \otimes \C{E}$, where $\C{Q}$ is as before but  
$\C{E}$ is a \emph{finite-dimensional} inner-product (Hilbert) space \citep{kreeSoize86}.
Such a situation arises when the parametric model is in
$\C{V} = \C{U} \otimes \C{E}$, where $\C{U}$ is an unspecified Hilbert
space as before, and one wants to see the `small' space
$\C{E}$ separately.  The parametric map can be defined as follows:
\begin{equation}   \label{eq:param-vect-r}
   \tns{r}:\C{M}\to\C{V} = \C{U} \otimes \C{E};\quad \tns{r}(\mu) = \sum_k r_k(\mu) \vek{r}_k,
\end{equation}
where as before $r_k(\mu) \in \C{U}$ and the $\vek{r}_k \in \C{E}$.  Typically the
index $k$ will range over the finite dimension of $\C{E}$, and the $\{\vek{r}_k\}_k$
are a suitable basis.  An example of such a situation is a vector field
over some manifold.  If at each point the vector is in the finite-dimensional space $\C{E}$,
e.g.\ the tangent space of the manifold, we model this by a tensor product
$\C{U} \otimes \C{E}$, where $\C{U}$ is some Hilbert space of scalar valued functions
on the manifold.

The `correlation' can now be given by a bilinear form.
The densely defined map $C_{\C{E}}$ in $\C{V}=\C{U}\otimes\C{E}$ is defined
on elementary tensors $\tns{u} = u\otimes\vek{u}, \tns{v}=v\otimes\vek{v}
 \in \C{V}=\C{U}\otimes\C{E}$ as
\begin{equation}   \label{eq:vect-corr}
   \bkt{C_{\C{E}}\tns{u}}{\tns{v}}_{\C{U}} :=   \sum_{k,j}  \bkt{R_k(u)}{R_j(v)}_{\C{Q}}
    \,(\vek{u}^{\ops{T}}\vek{r}_k)\, (\vek{r}_j^{\ops{T}}\vek{v})
\end{equation}
and extended by linearity, where each $R_k:\C{U}\to\C{Q}$ is the map associated to $r_k(\mu)$
as before for just a single map $r(\mu)$.
It may be called the 
\emph{`vector correlation'}.  By construction it is self-adjoint and positive.
The corresponding kernel with values in $\C{E}\otimes\C{E}$ for the
eigenvalue problem on $\C{W}=\C{Q}\otimes\C{E}$ is
\begin{equation}   \label{eq:vect-kernel}
   \vek{\vkappa}_{\C{E}}(\mu_1,\mu_2) = \sum_{k,j} \bkt{r_k(\mu_1)}{r_j(\mu_2)}_{\C{V}} 
   \; \vek{r}_k\otimes\vek{r}_j  .
\end{equation}

While the situation just described occurs often when $\tns{r}(\mu)$ is for
example the state of a system, a case which looks formally similar but 
allows an alternative approach happens
when the vector space $\C{E}$ consist of tensors of \emph{even} degree,
hence $\C{E} = \C{F}\otimes\C{F}$ for some space of tensors $\C{F}$ of half the degree.
An example of such a situation is the stress or strain field of a continuum
mechanics problem, where again $\C{U}$ could be a space of scalar spatial
functions, and $\C{E}$ the space of symmetric 2nd degree tensors.  Another
example is the specification of a conductivity tensor field for a heat conduction
problem with $\C{E}$ again the space of symmetric 2nd degree tensors, or
the specification of an elasticity tensor field with $\C{E}$
the space of symmetric 4th degree tensors.

Such a tensor of even degree can always be thought of as a linear map from a
space of tensors of half that degree into itself.  Being a linear map,
it can be represented as a \emph{matrix} $\vek{A}\in\D{R}^{n\times n}$, the
case we shall look at here.  The size of the matrix is equal to the
dimension of the space $\C{F}$, i.e.\ $n=\dim\,\C{F}$.
So the space $\C{E}$ can be thought of as a space of matrices.

Often these linear maps / matrices possess some
additional properties, like being symmetric positive definite --- as in the examples
of the conductivity or elasticity tensor --- or e.g.\ orthogonal.
One has to realise that the
representation methods which have been investigated here are \emph{linear} methods,
i.e.\ they work best when the representation is in a 
\emph{linear} manifold, essentially free from nonlinear constraints.
The two examples of positive definite or orthogonal tensor fields have such
nonlinear constraints.  We discuss two possible approaches:

First, assume that $\vek{A}$ has to be orthogonal.  It then satisfies
$\vek{A}^{\ops{T}}\vek{A} =\vek{I}=\vek{A}\vek{A}^{\ops{T}}$, a nonlinear constraint.
But the
orthogonal matrices $\ops{O}(n)$ and its sub-group of special orthogonal matrices
$\ops{SO}(n)$ form compact Lie groups.  One then may work in
their Lie algebra $\F{o}(n) = \F{so}(n)$, the \emph{skew} symmetric matrices,
the tangent space at the group identity $\vek{I}$.  This is a \emph{free} linear space.
An element $\vek{A}\in \ops{SO}(n)$ can be expressed with the exponential map
$\vek{A} = \exp(\vek{S})$ with $\vek{S}\in\F{so}(n)=\C{E}$.
Using the exponential map from the Lie algebra to its corresponding
Lie group one only has to deal with representations in the Lie algebra.

As a second example, assume that the matrix $\vek{A}\in\ops{Sym}^+(n)$ 
has to be symmetric positive definite (spd).  Then it can be
factored as $\vek{A}=\vek{G}^{\ops{T}}\vek{G}$ with invertible $\vek{G}\in\ops{GL}(n)$.
Both of these are nonlinear constraints.  The spd matrices $\ops{Sym}^+(n)$
are not a linear space, but
geometrically a salient open cone and a Riemannian manifold in the space of
all symmetric matrices $\F{sym}(n)$.
The manifold $\ops{Sym}^+(n)$ can be made into a Lie group in different ways.
Here it is important to observe that any $A\in\ops{Sym}^+(n)$ can be represented
again with the matrix exponential as $\vek{A}=\exp(\vek{H})$ with
$\vek{H}\in\F{sym}(n)=\C{E}$.  This is a good strategy
even for scalar fields, i.e.\ $n=1$.

Hence, in both cases, one may investigate the representation in some
\emph{linear sub-space}  $\F{g}\subseteq\C{F}\otimes\C{F}$, in the concrete matrix case
here in a sub-space  $\F{g}\subseteq\D{R}^{n \times n} = \F{gl}(n)$.
Such a parametric element may be represented first as $\vek{H}(\mu)\in\C{Q}\otimes\F{g}$
and then exponentiated:
\begin{equation}   \label{eq:exp-rep-spd}
   \vek{H}(\mu) = \sum_k \vsigma_k(\mu) \vek{H}_k,\qquad
   \vek{H}(\mu) \mapsto \exp(\vek{H}(\mu)) = \vek{A}(\mu) .
\end{equation}
Hence we now concentrate on representing $\vek{H}(\mu)$.
The parametric map would be written analogous to \refeq{eq:param-vect-r} as
\begin{equation}   \label{eq:tens-r-param}
   \tns{R}(\mu) = \sum_k r_k(\mu)\otimes\vek{R}_k \in \C{U}\otimes \C{E},\quad
   \text{ with } \quad \vek{R}_k \in \F{g}.
\end{equation}
The correlation analogous to \refeq{eq:vect-corr} may now be defined
via a bilinear form on elementary tensors as
a densely defined map $C_{\C{E}}$ in $\C{W}=\C{U}\otimes\C{F} = \C{U}\otimes\D{R}^n$ 
 --- observe, \emph{not} $\C{U}\otimes\C{E}=\C{U}\otimes\F{g}$ ---
and extended by linearity:
\begin{multline}   \label{eq:tens-corr}
\forall (\tns{u} = u\otimes\vek{v}), (\tns{v}=v\otimes\vek{v}) \in \C{W}=\C{U}\otimes\C{F}:\\
   \bkt{C_{\C{F}}\tns{u}}{\tns{v}}_{\C{U}} :=  \sum_{k,j}  \bkt{R_k(u)}{R_j(v)}_{\C{Q}}
    \,(\vek{R}_k\vek{u})^{\ops{T}}(\vek{R}_j \vek{v}).
\end{multline}
The kernel corresponding to \refeq{eq:vect-kernel} is again matrix valued,
\begin{equation}   \label{eq:tens-kernel}
   \vek{\vkappa}_{\C{F}}(\mu_1,\mu_2) = \sum_{k,j} \bkt{r_k(\mu_1)}{r_j(\mu_2)}_{\C{U}} 
   \; \vek{R}_k^{\ops{T}} \vek{R}_j  ,
\end{equation}
defining an eigenproblem in $\C{Q}\otimes\C{F}$.

For coupled systems, the approach has some similarity with the vector case $\C{U}\otimes\C{E}$
above.  The main characteristic of a coupled system which we want to preserve is that
the space state may be written as $\C{U}=\C{U}_1\times\C{U}_2$ with the the natural inner
product $\bkt{u}{v}_{\C{U}}=\bkt{u_1}{v_1}_{\C{U}_1}+\bkt{u_2}{v_2}_{\C{U}_2}$, where
$\vek{u}=(u_1,u_2), \vek{v}=(v_1,v_2) \in \C{U}$.  This is for two coupled systems, labelled as
`1' and `2'.  The parametric map  is
\begin{equation}   \label{eq:param-coup-r}
   \vek{r}:\C{M}\to\C{U} = \C{U}_1 \times \C{U}_2;\quad 
     \vek{r}(\mu) = (r_1(\mu), r_2(\mu)).
\end{equation}
The associated linear map is
\begin{equation}   \label{eq:lin_map-coup}
   \vek{R}:\C{U}\to\C{Q}^2 = \C{Q} \times \C{Q};\quad 
     (\vek{R}(\vek{u}))(\mu) = (\bkt{u_1}{r_1(\mu)}_{\C{U}_1}, \bkt{u_2}{r_2(\mu)}_{\C{U}_2}).
\end{equation}
As before, these $\D{R}^2$ valued functions on $\C{M}$ are like two problem-adapted
co-ordinate systems on the joint parameter set, one for each sub-system.
From this one obtains the `coupling correlation', again defined through a bilinear form
\begin{equation}   \label{eq:coup-corr}
   \bkt{\vek{C}_{c} \vek{u}}{\vek{v}}_{\C{U}} :=   \sum_{j=1}^2  \bkt{R_j(u_j)}{R_j(v_j)}_{\C{Q}} .
\end{equation}
The kernel is then a $2 \times 2$ matrix valued function in an integral operator
on $\C{W}=\C{Q}\times\C{Q}$:
\begin{equation}   \label{eq:coup-kernel}
   \vek{\vkappa}_{c}(\mu_1,\mu_2) =  
   \diag (\bkt{r_k(\mu_1)}{r_k(\mu_2)}_{\C{U}_k}) .
\end{equation}

Often there is a bit more structure one wants to preserve, namely that
$\C{M}=\C{M}_1\times\C{M}_2$, and the parameter set $\C{M}_1$ is for the sub-system `1',
and the set $\C{M}_2$ is for sub-system `2'.  We also assume that not only 
$\C{U}=\C{U}_1\times\C{U}_2$, but also $\C{Q}=\C{Q}_1\times\C{Q}_2$, where the scalar
functions in $\C{Q}_1$ depend only on $\C{M}_1$, and similarly for subsystem `2'.
The parametric map is hence
\begin{equation}   \label{eq:param-coup-12}
   \vek{r}:\C{M}=\C{M}_1\times\C{M}_2\to\C{U} = \C{U}_1 \times \C{U}_2;\quad 
     \vek{r}((\mu_1,\mu_2)) = (r_1(\mu_1), r_2(\mu_2)) ,
\end{equation}
with the associated linear map
\begin{equation}   \label{eq:lin_map-coup-12}
   \vek{R}:\C{U}\to\C{Q} = \C{Q}_1 \times \C{Q}_2;\quad 
     (\vek{R}(\vek{u}))(\mu) = (\bkt{u_1}{r_1(\mu_1)}_{\C{U}_1}, \bkt{u_2}{r_2(\mu_2)}_{\C{U}_2}).
\end{equation}
The correlation may be defined as before in \refeq{eq:coup-corr}, and also the
kernel on $\C{Q}=\C{Q}_1\times\C{Q}_2$ is as in \refeq{eq:coup-kernel},
but now the first diagonal entry is a function on $\C{M}_1 \times \C{M}_1$ only,
and analogous for the second diagonal entry.

%
%
%
%
%
%
%
%
%
%
%
%


%

\section{Conclusion} \label{S:concl}
Parametric mappings $r:\C{M}\to\C{U}$ have been analysed with in a  variety of
settings via the associated linear map $R:\C{U}\to\C{Q}\subseteq\D{R}^{\C{M}}$,
enabling the linear analysis.  The RKHS setting allows a first representation, and
essentially reproduces everything in $\C{U}$ in the function space $\C{R}$.
The choice of another inner product and corresponding Hilbert space $\C{Q}$
leads to measures of importance in  $\C{M}$, or, more precisely, in $\D{R}^{\C{M}}$.

It is shown that each separated representation defines an associated linear map,
and that conversely under some more restrictive conditions, the normally more general
notion of an associated linear map defines a representation.

Several refinements are presented to represent some additional structure in the
linear map.  One such structure is the information of dealing with a coupled system.
This can be reflected in the structure of the associated linear map.

\ignore{           
\subsection{Test of fonts} \label{SS:testf}
\paragraph{Slanted:} --- serifs
\[ a, \alpha, A, \Phi, \phi, \vphi, \qquad 
     \vek{a}, \vek{\alpha}, \vek{A}, \vek{\Phi}, \vek{\phi}, \vek{\vphi} \]

\paragraph{Slanted:} --- sans serif
\[ \tns{a}, \tns{\alpha}, \tns{A}, \tns{\Phi}, \tns{\phi}, \tns{\vphi}, \qquad
  \tnb{a}, \tnb{\alpha}, \tnb{A}, \tnb{\Phi}, \tnb{\phi}, \tnb{\vphi} \]
   
\paragraph{Other:}  --- no small letters
\[ \C{E}, \C{Q}, \C{R}; \D{E}, \D{Q}, \D{R}; \E{E}, \E{Q}, \E{R}\]

\paragraph{Fraktur:}
\[ \F{e}, \F{q}, \F{r}; \F{E}, \F{Q}, \F{R} \]

\paragraph{Upright:} --- serifs, no small greek letters
\[ \mrm{a}, \mrm{A}, \mrm{\Phi},\qquad \mat{a}, \mat{A}, \mat{\Phi} \]

\paragraph{Upright:} --- sans serif, no small greek letters
\[ \ops{a}, \ops{A}, \ops{\Phi},\qquad \opb{a}, \opb{A}, \opb{\Phi} \]
}           

%
%
%
%
%
%
%
%
%
%
%





\bibliography{\thebib/jabbrevlong,\thebib/stochastic,\thebib/fuq-new,%
\thebib/fa,\thebib/mat_BU-1-S,\thebib/phys_D,\thebib/num,\thebib/highdim}


{ 
   \tiny
       \texttt{\RCSId} 
   }



\end{document}